\documentclass[12pt,twoside]{article}

\usepackage{graphicx,graphics, amsthm, amsfonts, amsbsy,amssymb, amsmath, cite, array,arydshln, hyperref, float,tikz}
\usepackage[utf8]{inputenc}

\let\OLDthebibliography\thebibliography
\renewcommand\thebibliography[1]{
	\OLDthebibliography{#1}
	\setlength{\parskip}{0pt}
	\setlength{\itemsep}{0pt plus 0.3ex}
}

\everymath{\displaystyle}

\usepackage[margin=1in]{geometry}

\allowdisplaybreaks

\newtheorem{theorem}{Theorem}[section]

\newtheorem{proposition}[theorem]{Proposition}
\newtheorem{note}{Note}

\renewenvironment{proof}{{\bfseries Proof.}}{\qed}

\begin{document}
\title{A note on the Independent domination polynomial of zero divisor graph of rings}
\author{Bilal Ahmad Rather\footnote{Orcid:0000-0003-1381-0291}\\
	{\em Mathematical Sciences Department, College of Science,}\\
	\emph{United Arab Emirates University, Al Ain 15551, Abu Dhabi, UAE}\\
	bilalahmadrr@gmail.com
				}
\date{}

\pagestyle{myheadings} \markboth{Bilal Ahmad Rather}{A note on the Independent domination polynomial of zero divisor graph of rings}
\maketitle

\begin{abstract}
In this note we consider the independent domination polynomial problem along with their unimodal and log-concave properties which were earlier studied by G\"ursoy, \"Ulker and G\"ursoy (Soft Comp. 2022). We show that the independent domination polynomial of zero divisor graphs of $\mathbb{Z}_{n}$ for $n\in \{ pq, p^{2}q, pqr, p^{\alpha}\}$ where $p,q,r$ are primes with $2<p<q<r$ are not unimodal thereby contradicting the main result of G\"ursoy, \"Ulker and G\"ursoy \cite{gursoy}. Besides the authors show that the zero of the independent domination polynomial of these graphs have only real zero and used concept of Newton's inequalities to establish the log-concave property for the afore said polynomials. We show that these polynomials have complex zeros and the technique of  Newton's inequalities are not applicable. Finally, by definition of log-concave, we prove that these polynomials are  log-concave and fix the flaws in Theorem 10 of G\"ursoy, \"Ulker and G\"ursoy \cite{gursoy}. 
\end{abstract} 
\vskip 3mm



\section{Introduction}
\paragraph{}
All graphs considered in this note  are finite, simple and undirected graphs.  A graph is usually symbolized by $G=G(V(G),E(G))$ with its vertex set $V(G)$ and edge set $E(G)$. The numbers $n = |V(G)| $ is order and $ m = |E(G)| $ size of $ G $  An edge among vertices $ u $ and $ v $ is denoted by $ u\sim v $. 
The \textit{degree} $ d_{v_{i}}(G) $ of a vertex $ v_{i}\in V(G) $ is the number of vertices incident on it. 
The \textit{union} of two graphs $ G_{1}(V_{1}(G_{1}), E_{1}(G_{1})) $ and $ G_{2}(V_{2}(G_{2}), E_{2}(G_{2})) $, denoted by  $ G_{1}\cup G_{2}$,  defined as a graph with vertex set $ V_{1}(G_{1})\cup V_{2}(G_{2}) $ and edge set $ E_{1}\cup E_{2}. $ The \textit{join} of $ G_{1}(V_{1}(G_{1}), E_{1}(G_{1})) $ and $ G_{2}(V_{2}(G_{2}), E_{2}(G_{2})) $, denoted by $ G_{1}\vee G_{2} $,   is  a graph with vertex set $V_{1}(G_{1})\cup V_{2}(G_{2})$ and edge set $E(G_{1})\cup E(G_{2})\cup\big\{u,v~|~ u\in V(G_{1}), v\in V(G_{2}) \big\}.$

A set $\emptyset\neq  S\subseteq V(G) $ is called a \emph{dominating} set if each vertex in $ V\setminus S $ is adjacent to at least one  vertex in $ S. $ The minimum order (cardinality) among all such dominating sets of $ G $ is called the \emph{domination number}  $ \gamma(G) $ of $ G. $ The theory of domination of graphs is well studied, see textbook \cite{haynes}. A set $S\subseteq V(G)$ in a graph $G$ is called an \emph{independent} set if vertices of $S$ are pairwise non-adjacent. The cardinality of the maximum independent set is called the \emph{independence number} $ \alpha(G)$  of $ G. $
A subset $D$ of $V(G)$ is independent dominating set of $ G $ which is both dominating and independent in $ G $. The independent domination number $ \gamma_{i}(G) $ is the minimum size of all independent dominating sets of $ G. $ The invariants $ \gamma, \alpha $ and $ \gamma_{i} $ of $ G $ is related by $ \gamma(G)\leq \gamma_{i}(G) \leq \alpha(G)$ (see, \cite{haynes}). The independent set problem in a graph is strongly NP-hard problem while the dominating set problem of a graph is NP-complete problem. These problems are well studied both in theoretical computer science and mathematics. 

A polynomial $ p(x)=\sum_{i=0}^{b} a_{i}x^{i} $ is called \textit{unimodal} if  coefficients $a_{i}$'s form a unimodal sequence, that is,  there exists a positive integer  $ p~(0\leq p\leq n), $ known as the \emph{mode}, such that  $ a_{0}\leq a_{1}\leq \dots \leq a_{p} \geq a_{p+1}\geq \dots \geq a_{b} $. Equivalently the coefficients of $p(x)$ increase to some stage and then decrease from thereafter.
The polynomial $ p(x) $ is \emph{log-concave} (logarithmically concave) if 
\begin{equation}\label{log con 1}
	a_{j}^{2}\geq a_{j-1}a_{j+1},~  \text{for all}~ 1\leq j\leq b-1. 
\end{equation} 

If $a_{i}$'s are non negative and all the zeros of $p(x)$ are real. Then the basic approach for unimodal and log-concave  is the  Newton's inequalities \cite{hardy} stated as:
\begin{equation}\label{log con 2}
	a_{i}^{2}\geq a_{i-1}a_{i+1}\left(1+\frac{1}{i}\right)\left(1+\frac{1}{b-i}\right),~\text{for}~ i=1,\dots,b-1. 
\end{equation}

The condition for log-concave given in \eqref{log con 2} is stronger than the one given in \eqref{log con 1}, see \cite{stanley}. Also log-concave sequence of positive terms is unimodal \cite{stanley}. However, if $a_{i}=0$ implies that either $a_{1}=\dots=a_{i-1}=0$ or $a_{i+1}=\dots=a_{n}=0,$ (no internal zeros), then a non-negative log-concave sequence is unimodal \cite{stanley}.
 Log-concavity of sequences is related to the surface embedding of graphs in topological graph theory, like the log-concavity of genus polynomials of  graphs \cite{gross}.

Clearly, the polynomial  $ 1 + 2x + 3x^{2} + 7x^{3} + 6x^{4} + 3x^{5} $ is unimodal whereas $ 21 + 7x + 19x^{2} + 6x^{3} + 4x^{4} $ is not unimodal, since the coefficients decrease, then increase and then decrease. An effective information about the increase or decrease of coefficients can measure the unimodal property. As it is clear that unimodal sequence either increases or decreases, or increase and then decrease. The  number of changes of directions (increasing or decreasing) of $p(x)$ is defined as the of oscillations $ \eta(p(x)) $ of $ p(x) $. Obviously,  for the unimodal polynomial the oscillations of  $p(x)$ must be at most one. 
For polynomial $ p(x) = 8 + 6x + 160x^{2} + 28x^{3} +16x^{4} +11x^{5} $, then $ \eta(p(x))=2 $, since change of decreasing directions are two. Sometimes it is equally good to identify the oscillations for non unimodal polynomials.

Identifying unimodal or log-concave property (or both) of a given polynomial  is a very well studied non-trivial  problem. Several open problems and conjectures exist in literature related to these entities for different types of polynomials. 
These properties of polynomials are mostly studied for graph polynomial, like independent polynomial, dominating polynomial, matching polynomial, clique polynomial and several other polynomials.

Let $ d_{k}(G,k) $ be the number of independent dominating sets of order $k$ in $ G $. The independent domination polynomial of $ G $ is defined as
\[ D_{i}(G,x)=\sum_{k=\gamma_{i}(G)}^{\alpha(G)}d_{i}(G,k)x^{k}. \]
A root of the equation $ D_{i}(G,x)=0 $ is called as the independent domination root of $ G. $

 The independent domination polynomial $ D_{i}(G,x) $ is a generating function of number of the independent dominating sets of certain cardinalities of $ G $. 
The independent domination polynomials and their roots, unimodal and log-concave property have attracted many researchers, see \cite{anwar,lonzaga, dod, goddard}. Jahari and Alikhani \cite{jahari} obtained the independent domination polynomials of generalized compound graphs and constructed graphs whose independent domination polynomials have real zeros.  More about independent domination polynomials and other polynomials can be seen in \cite{alikhani, alikhani1, alikhani2, levitdiscrete, levitsurvey, anthony}. Recently, G\"ursoy, \"Ulker and G\"ursoy \cite{gursoy} presented the results related to the independent domination polynomial of zero-divisors graphs associated to commutative rings. In particular, they showed  that the independent domination polynomial of zero divisor graph of $\mathbb{Z}_{n}$  is unimodal and log-concave. We will consider this study in the present note and modify the existing results of G\"ursoy, \"Ulker and G\"ursoy \cite{gursoy} and fix the errors in the published article.

In Section \ref{section 2}, we give basic of zero divisor graphs of $\mathbb{Z}_{n}$ and recall the results of G\"ursoy, \"Ulker and G\"ursoy \cite{gursoy} along with the properties of log-concave and unimodal. We give several examples having complex independent domination zeros countering a result of \cite{gursoy}. We discuss the unimodal and log-concave property of $\Gamma(\mathbb{Z}_{n})$ for $n\in \{p^{2}, 2p, pq, p^{2}q, pqr\}$ where $p,q,r$ are primes in Propositions \ref{unimodal prime p}, \ref{unimodal prime pq} and \ref{unimodal prime pq^{2}}. Theorems \ref{thm zpalpha} gives the independent domination polynomial in simplified form and Theorem \ref{thm zpalpha restate} established the non unimodal and log-concave property of $\Gamma(\mathbb{Z}_{p^{\alpha}})$ where $p>2$ is prime and $\alpha$ is  a positive integer.


\section{Independent domination polynomial of zero divisor graphs of commutative rings}\label{section 2}
\paragraph{}

For a given commutative ring  $ R $  with non-zero identity, $ 0\neq a\in R $ is a zero divisor of $ R $ if there is $ 0\neq b\in R $ such that $ a\cdot b=0. $ Beck in 1988 \cite{ib} put forward the concept of zero divisor graphs $\Gamma(R)$ to study the colouring of rings. He included identity $0 $ in the vertex set of zero divisor graphs and defined edges among zero divisor if and if their product is zero, thereby we see that $0$ is adjacent to all vertices of such a graphs.
Anderson and  Livingston in 1999 \cite{al} modified the definition Beck's zero divisor graphs and excluding $ 0 $ of ring in zero divisor graph and defined edges between two non-zero zero divisors if and only if their product is zero. Several interesting results related to these graphs and their underlying rings were published in the past years, mostly recent notably works are: \ Chattopadhyay  Panigrahi \cite{sc} obtained the graphical structure of $\Gamma(\mathbb{Z}_{n})$ and showed that it is join of certain complete and null graphs (non empty graphs without any edges). Anderson and Weber \cite{anderson2019} determine all zero-divisor graphs with at most $14$ vertices.
First spectral analysis of zero divisor graphs were carried in  \cite{my}. Complete spectral and structural analysis were given in \cite{bilalmathematics}. G\"ursoy, \"Ulker and G\"ursoy  \cite{gursoy} introduced the concept of independent domination polynomials for zero divisor graphs for $\mathbb{Z}_{n}$ for certain values of $n.$ We will carry the work G\"ursoy, \"Ulker and G\"ursoy  \cite{gursoy} forward in this note. First we understand the structure of  $\Gamma(\mathbb{Z}_{n})$.

Consider the sets as given below:
\[ V_{e_{i}}= \{ x\in \mathbb{Z}_{n} : (x,n)=e_{i} \}, \quad \text{for}\quad  1\leq i \leq t, \]
where $t\notin \{1,n\} $ is a divisor of $n$ and $ (x,n) $ denotes the greatest common divisor of $ x $ and $ n. $ It is clear that $ A_{e_{i}} \cap A_{e_{j}}=\emptyset$, for $ i\neq j $. Thus, it implies that $ A_{e_{1}}, A_{e_{2}}, \dots, A_{e_{t}} $ are pairwise disjoint and partitions the vertex set of $ \Gamma(\mathbb{Z}_{n}) $ as
\[
V(\Gamma(\mathbb{Z}_{n}))=A_{e_{1}}\cup A_{e_{2}}\cup \dots \cup A_{e_{t}},
\]
and vertices of $ A_{e_{i}} $ are adjacent to vertices of $ A_{e_{j}}  $ in if and only if $ e_{i}e_{j}=kn, $ where $k$ is scaler (see \cite{sc}). Also, 
 $|A_{e_{i}}|=\phi\left( \frac{n}{e_{i}}\right)$, for $ 1\leq i \leq t $ (see \cite{my}). Besides the induced subgraphs of $ A_{e_{i}} $ is either  either clique or its complement. More precisely $ \Gamma(E_{d_i}) $ is $ K_{\phi\left (\frac{n}{d_{i}}\right )} $ is $ e_{i}^{2} $ divides $ n $ otherwise it is $\overline{K}_{\phi\left (\frac{n}{d_{i}}\right )}$ (see \cite{sc}). Next, we state a results related to the unimodal and log-concave property of the independent domination polynomial of $\Gamma(\mathbb{Z}_{n})$ for several values of $n,$ proved as a main and interesting result in G\"ursoy, \"Ulker and G\"ursoy  \cite{gursoy}. We restate it and avoid the ambiguity in the stated version of \cite{gursoy}.

\begin{theorem}[\cite{gursoy}, Theorem 10]\label{gursoy}
	The independent domination polynomials of  zero divisor graphs of rings $\mathbb{Z}_{pq} , \mathbb{Z}_{p^{2}q} $ and $ \mathbb{Z}_{pqr}$ is unimodal and log-concave for prime numbers $p, q, r ,$ with $ p > q > r$. The independent domination polynomial of zero divisor graph of $\mathbb{Z}_{p^{\alpha}} $ is log-concave and unimodal for prime $p>2$, where  $\alpha$ is a positive integer.
\end{theorem}
The proof of Theorem 10 from G\"ursoy, \"Ulker and G\"ursoy  \cite{gursoy} reads the following:\vskip 1.5mm

\textit{The independent domination polynomials of  $\mathbb{Z}_{pq} , \mathbb{Z}_{p^{2}q} , \mathbb{Z}_{pqr}$ and $\mathbb{Z}_{p^{\alpha}}$ are unimodal since the sequence of its coefficients $\alpha_{0}, \alpha_{1},\dots, \alpha_{m}$ are a sequence of non-negative numbers, and these polynomials have only real zeros. Then, we have 
$$\alpha^{2}_{t}\geq \alpha_{t-1}\alpha_{t+1}\left(1 + \frac{1}{t}\right)\left(1 + \frac{1}{m-t}\right), t = 1, 2,\dots ,m -1$$ from Newton’s inequalities. Therefore, the sequence is log-concave and unimodal.}

\begin{note}
We note two important points from the above result and its proof.
\end{note}
\begin{enumerate}
	\item The zero of the independent domination polynomials of the zero divisor graphs of $\mathbb{Z}_{pq} , \mathbb{Z}_{p^{2}q} , \mathbb{Z}_{pqr}, $ and $ \mathbb{Z}_{p^{\alpha}}$  are not always real and they have complex zeros as well, where  $p, q, r ,$ are primes with $p > 2$ and $\alpha$ is a positive integer. We will show it in next couple of examples along with their graphical representation on complex plane.
	\item Since the zeros of the independent domination polynomials of the zero divisor graphs of $\mathbb{Z}_{pq} , \mathbb{Z}_{p^{2}q} , \mathbb{Z}_{pqr}$ and $\mathbb{Z}_{p^{\alpha}}$ are not always real, so Newton’s inequalities are not applicable and log-concave cannot be established with this procedure used in the proof pf Theorem 10 of \cite{gursoy}.
\end{enumerate}

Next, we give a sequence of examples of the independent domination polynomials of the zero divisor graphs of the  families of the commutative rings mentioned in Theorem \ref{gursoy} having complex zeros as well and are not unimodal. We consider the some examples from \cite{gursoy}. Their graphs are given in \cite{gursoy} and graphs of other examples can be drawn in a similar fashion. 

\begin{enumerate}
	\item For $G\cong \Gamma(\mathbb{Z}_{pq})$ with $p=3$ and $q=5,$  the independent domination polynomial is $D_{i}(G,x)=x^{2}+x^{4}$ (see Theorem 4, \cite{gursoy}) and its zeros are $0$ and $\pm i$. The graphical representation of  the zeros of $D_{i}(G,x)$ is easy to visualise. 

\item For $G\cong \Gamma(\mathbb{Z}_{p^{2}q})$ with $p=5^{2}$ and $q=3,$  the independent domination polynomial is $D_{i}(G,x)=x^{10}+4x^{21}+x^{28}$ (see Example 1, \cite{gursoy}) and its set of zeros are 
\begin{footnotesize}
	\begin{align*}
	\Big\{&0, 1.21378, -0.891122, -0.764225-0.955193 i, -0.764225+0.955193 i, -0.739951-0.469162 i,\\
	& -0.739951+0.469162 i, -0.358291-0.803288 i, -0.358291+0.803288 i, 0.120955\, +0.88068 i,\\
	& 0.120955\, -0.88068 i, 0.274622\, -1.18462 i, 0.274622\, +1.18462 i, 0.571168\, -0.662058 i,\\
	& 0.571168\, +0.662058 i, 0.850704\, +0.241222 i, 0.850704\, -0.241222 i, 1.09747\, -0.533693 i,\\
	& 1.09747\, +0.533693 i\Big\}.
\end{align*}
\end{footnotesize}
 
The graphical representation of the zeros of $D_{i}(G,x)$ is shown in Figure  \ref{Fig zeros pq}, where green dots represent the zeros. 
\begin{figure}[H]
	\centering{\scalebox{.285}{\includegraphics{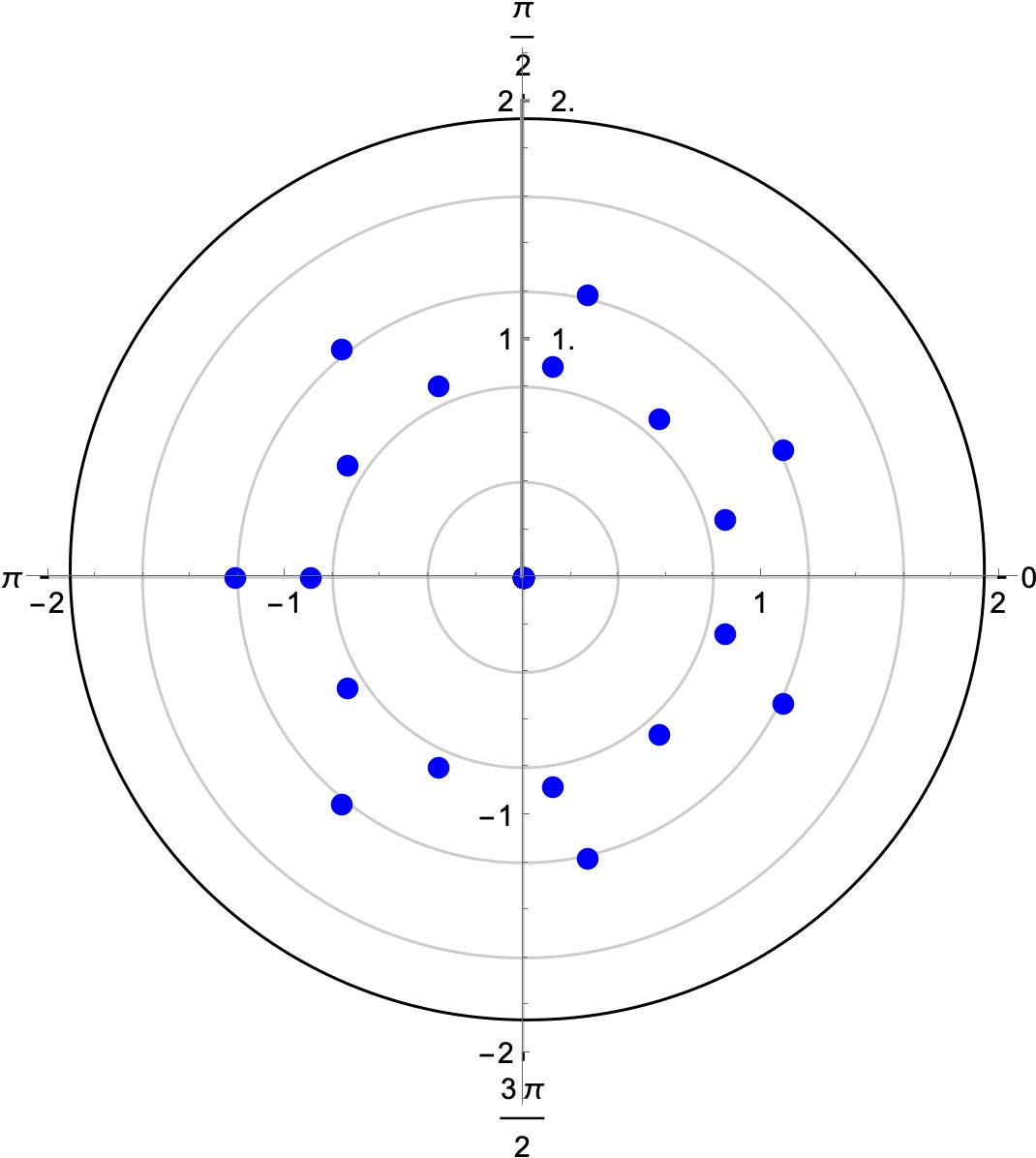}}\qquad  \scalebox{.285}{\includegraphics{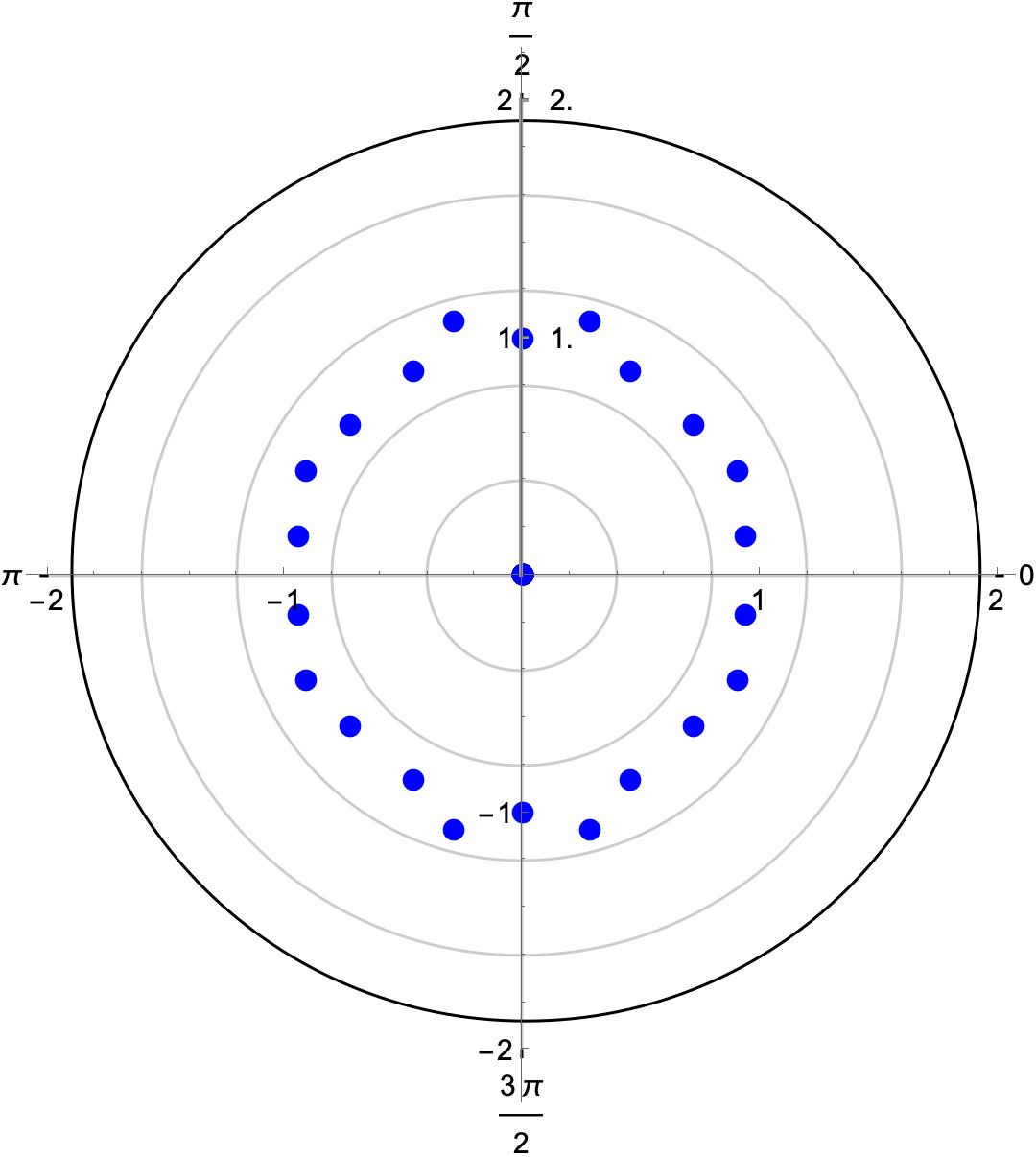}}}\\
	\qquad Zeros of $D_{i}(\Gamma(\mathbb{Z}_{75}), x).$\qquad\qquad  Zeros of $D_{i}(\Gamma(\mathbb{Z}_{105}), x).$
	\caption{Pictorial representation of the zeros of  on a plane.}
	\label{Fig zeros pq}
\end{figure}

\item For $G\cong \Gamma(\mathbb{Z}_{pqr})$ with $p=3, q=5$ and $r=7,$  the independent domination polynomial is $D_{i}(G,x)=x^{22}+x^{36}+x^{42}+x^{44}$ (see Example 2, \cite{gursoy}) and its set of zeros are 
\begin{footnotesize}
	\begin{align*}
		\Big\{&0,-0.937799-0.164137 i, -0.937799+0.164137 i, -0.907901-0.436924 i, -0.907901+0.436924 i,\\
		& -0.721752-0.635391 i, -0.721752+0.635391 i, -0.455901-0.860636 i, -0.455901+0.860636 i,\\
		& -0.29053-1.07456 i, -0.29053+1.07456 i,  -1. i, 0.\, +1. i, 0.29053\, -1.07456 i, 0.29053\, +1.07456 i,\\
		& 0.455901\, -0.860636 i, 0.455901\, +0.860636 i, 0.721752\, -0.635391 i, 0.721752\, +0.635391 i, 0.907901\,\\
		& -0.436924 i, 0.907901\, +0.436924 i, 0.937799\, -0.164137 i, 0.937799\, +0.164137 i\Big\}.
	\end{align*}
\end{footnotesize}
The graphical representation of the zeros of $D_{i}(G,x)$ is shown in Figure  \ref{Fig zeros pq} (right), where green dots represent the zeros.

\item For $G\cong \Gamma(\mathbb{Z}_{p^{\alpha}})$ with $p=3$ and $\alpha=5,$  the independent domination polynomial is $D_{i}(G,x)=2x+6x^{55}+x^{72}$ (see Example 3, \cite{gursoy}) and the set of zeros of this polynomial is very long. So, we represent it by the graphical picture as shown in Figure \ref{Fig zeros palpha} 

\begin{figure}[H]
	\centering{\scalebox{.285}{\includegraphics{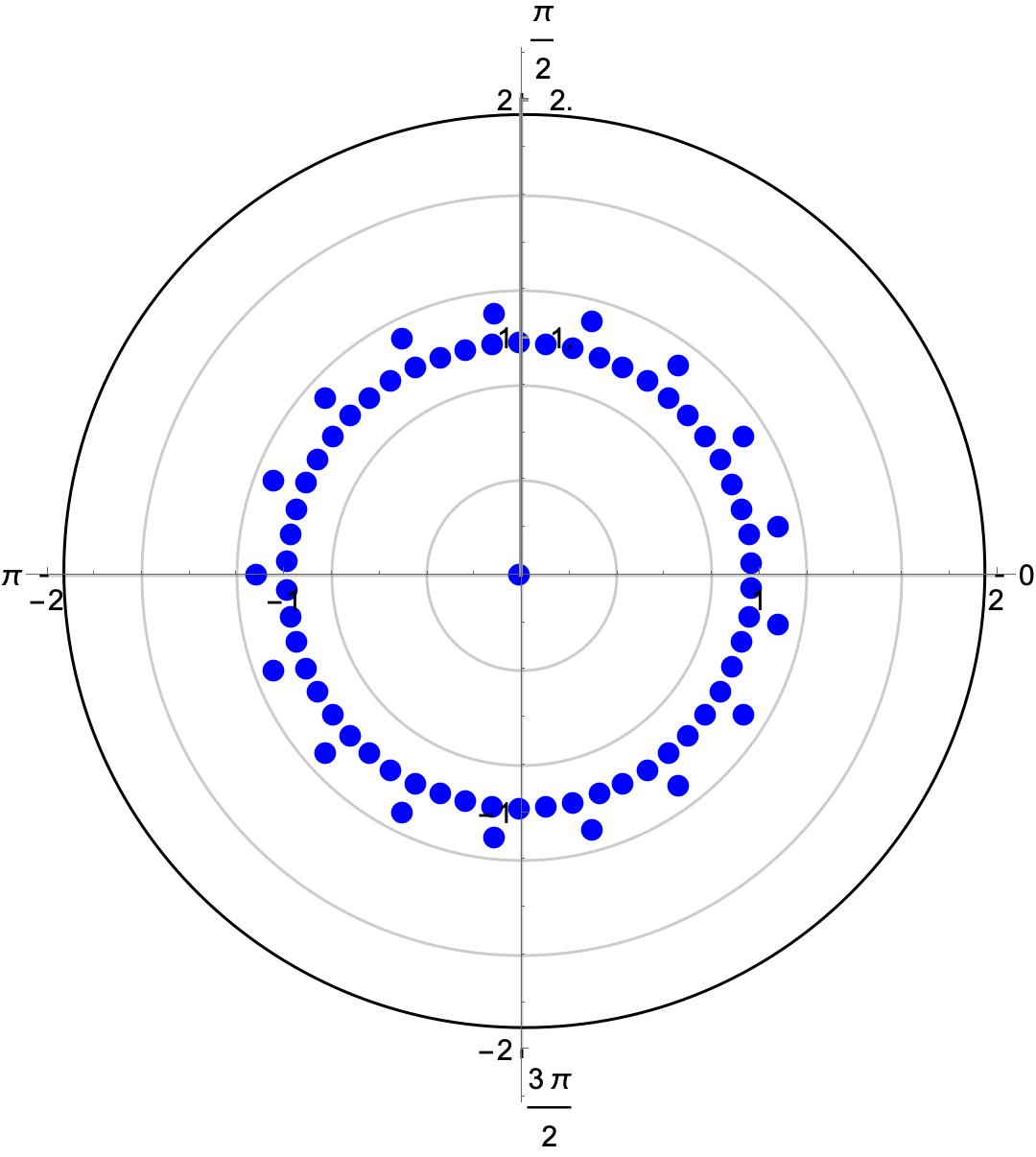}}\qquad  \scalebox{.285}{\includegraphics{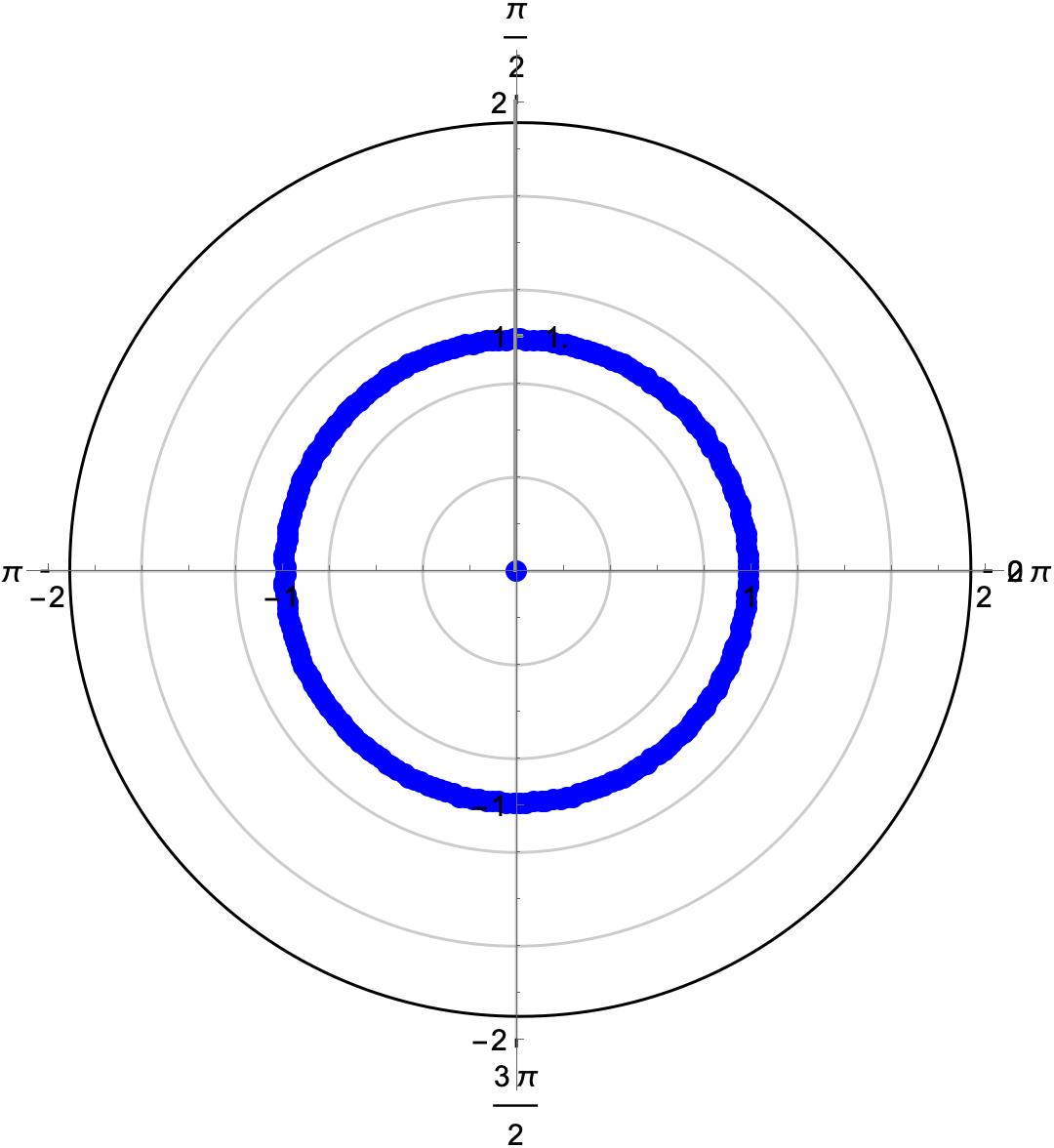}}}\\
	\qquad Zeros of $D_{i}(\Gamma(\mathbb{Z}_{243}), x).$\qquad\qquad  Zeros of $D_{i}(\Gamma(\mathbb{Z}_{729}), x).$
	\caption{Pictorial representation of the zeros of  on a plane.}
	\label{Fig zeros palpha}
\end{figure}

\item For $G\cong \Gamma(\mathbb{Z}_{p^{\alpha}})$ with $p=3$ and $\alpha=6,$  the independent domination polynomial is $D_{i}(G,x)=2x+6x^{163}+x^{216}+18x^{217}$ and the list of zeros of this polynomial is very long. So, we avoid listing them and present its graphical representation in Figure \ref{Fig zeros palpha} (right).
\end{enumerate}

Next, we will check the unimodal property of the independent domination polynomial of the zero divisor graphs for $\mathbb{Z}_{pq}, \mathbb{Z}_{p^{2}q}, \mathbb{Z}_{pqr}$ and $\mathbb{Z}_{p^{\alpha}}.$

We recall that the independent domination polynomial of $\Gamma(\mathbb{Z}_{p^{2}})$ is $(p-1)x,$ which is clearly unimodal as sequence of its coefficients form an increasing sequence of non-negative integers. From now onwards, whenever we consider $n=p^{\alpha}$, we assume $p>2$ to avoid $\Gamma(\mathbb{Z}_{p^{2}})$ case.
 Next for $D_{i}(\Gamma(\mathbb{Z}_{2p}),x)=x+x^{p-1}$ to be unimodal, we have $a_{0}=0\leq a_{1}=1$ which must be less than or equal to coefficient of $x^{p-1}$, which can happen if and only if $p-1=2$ as its coefficient is already one. Otherwise if $p>2,$ then some coefficients between the linear term and $(p-1)$-th must be zero. Which will imply that $\eta(D_{i}(\Gamma(\mathbb{Z}_{2p}),x))>2,$ thereby contradicting the requirement for unimodal property.  Also, by definition the log-concave property is satisfied trivially for $\Gamma(\mathbb{Z}_{p^{2}})$. However, for $D_{i}(\Gamma(\mathbb{Z}_{2p}),x)=x+x^{p-1},$ we see that $1=a_{1}^{2}\geq a_{0}a_{2}=0$ and for $i\geq 2$, we get $0=a_{i}^{2}\geq a_{i-1}a_{i+1} =0$, which are always true, and it follows that $D_{i}(\Gamma(\mathbb{Z}_{2p}),x)$ is log-concave. We make these observation precise in the following result.
\begin{proposition}\label{unimodal prime p}
	Let $p$ be a prime. Then the following hold.
	\begin{itemize}
		\item[\bf (i)] For prime $p$, the independent domination polynomial of $\Gamma(\mathbb{Z}_{p^{2}})$ is unimodal and log-concave. 
		\item[\bf(ii)] For prime $p>2$, the independent domination polynomial of $\Gamma(\mathbb{Z}_{2p})$ is always log-concave and is unimodal if and only if $p=3.$ 
	\end{itemize}
\end{proposition}

With a similar argument as in (ii) of Proposition \ref{unimodal prime p}, we have the following result.
\begin{proposition}\label{unimodal prime pq}
	The independent domination polynomial of $ \Gamma(\mathbb{Z}_{pq})$ with primes $p<q$ is log-concave but it is unimodal if and only if $p=2$ and $q=3.$
\end{proposition}

\begin{proposition}\label{unimodal prime pq^{2}}
	The independent domination polynomial of $ \Gamma(\mathbb{Z}_{p^{2}q})$ with primes $p>q$ is log-concave but not unimodal. 
\end{proposition}
\noindent\begin{proof}
	From Theorem 5 \cite{gursoy}, without evaluating Euler functions, the independent domination polynomial of $G\cong \Gamma(\mathbb{Z}_{p^{2}q})$  is
	\[ D_{i}(G,x)=x^{\phi(q)+\phi(pq)}+\phi(p)x^{\phi(p^{2})+1}+x^{\phi(p^{2})+\phi(pq)}. \]
	For $p=3$ and $q=2,$ the above polynomial becomes 
	\[D_{i}(G,x)=x^{3}+2x^{7}+x^{8},  \]
	and it is clear that $\eta(D_{i}(G,x))=2$, since the change of decreasing (increasing) directions of coefficients is two. So, it is not unimodal. Now, for other values of $p$ and $q$, if $\phi(q)+\phi(pq)<\phi(p^{2})+1$, then we obtain $p(q-1)<p^{2}-p<p^{2}-p+1$ and from it, we get $q<p $, which is true. Also,  $\phi(p^{2})+1\leq \phi(p^{2})+\phi(pq)$. Thus, the coefficients of $D_{i}(G,x)$ increase, then decrease and then again increase, then decrease. Thus, it implies that $\eta(D_{i}(G,x))>2$ and $D_{i}(G,x)$ is not unimodal. It is easy to see that the coefficients of $D_{i}(G,x)$ satisfies the conditions $a_{i}^{2}\geq a_{i-1}a_{i+1}$ for $1\leq i\leq \phi(p^{2})+\phi(pq)$ and it implies the log-concave property of $D_{i}(G,x). $
\end{proof}

We note that the number of oscillation in the independent domination polynomial of the above results is presciently three, since $a_{0}\leq \dots \leq a_{\phi(q)+\phi(pq)}\geq a_{\phi(q)+\phi(pq)+1}\geq \dots\geq a_{\phi(p^{2})} \leq a_{\phi(p^{2})+1}\geq a_{\phi(p^{2})+2}\dots\geq a_{\phi(p^{2})+\phi(pq)-1} \leq a_{\phi(p^{2})+\phi(pq)}. $

The following result gives the unimodal and log-concave property for the zero divisor graphs when $n$ is product of three primes and the result corrects Theorem 10 of \cite{gursoy} for $\Gamma(\mathbb{Z}_{pqr}) $.
\begin{proposition}
	The independent domination polynomial of $ \Gamma(\mathbb{Z}_{pqr}) $ with primes $p<q<r$ is never unimodal and it is log-concave if and only if $|\phi(pq)-\phi(r)|\neq 2.$
\end{proposition}
\noindent\begin{proof}
	Let $G\cong  \Gamma(\mathbb{Z}_{pqr})$ be the zero divisor graph of $\mathbb{Z}_{pqr}$, where we assume that $p<q<r$ without any loss of generality. The independent domination polynomial of $G$ is (see Theorem 6, \cite{gursoy})
	\begin{equation}\label{eq pqr}
		D_{i(G,x)}=x^{\phi(pr)+\phi(pq)+\phi(p)}+x^{\phi(qr)+\phi(pq)+\phi(q)}+x^{\phi(qr)+\phi(pr)+\phi(r)}+x^{\phi(qr)+\phi(pr)+\phi(pq)}.
	\end{equation}
	If $p\neq 2,$ then all the exponents of the polynomial given in \eqref{eq pqr} are Euler function, which are even and so is their sum. Thus, $\eta(D_{i}(G,x))=4$ as coefficients increase four times and decrease three times. So, $D_{i}(G,x)$ is not unimodal. If $p=2$, then exponent $\phi(pr)+\phi(pq)+\phi(p)$ in \eqref{eq pqr} is odd and other coefficients are even. In this case $\eta(D_{i}(G,x))$ is at least three and it implies that $D_{i}(G,x)$ is not unimodal. However depending on $p<q<r$, it is clear that either $ \phi(qr)+\phi(pr)+\phi(r)$ or $\phi(qr)+\phi(pr)+\phi(pq)$ is the largest exponent. Also, if $\phi(pr)+\phi(pq)+\phi(p)<\phi(qr)+\phi(pq)+\phi(q)$ with $p=2$, then $\phi(q)>\phi(p)=1$ and $\phi(qr)>\phi(2r)=\phi(pr)=\phi(r)$. It implies that there is difference of more two between $\phi(qr)+\phi(pq)+\phi(q)$ and $\phi(pr)+\phi(pq)+\phi(p).$ Thus, in this case we obtain $\eta(D_{i}(G,x))=4$ and polynomial is not unimodal. Next, we turn our attention toward the log-concave property of $D_{i}(G,x)$. It is clear that the minimum value of $r$ is $5$, so $\phi(r)-\phi(q)\geq 2$ and $\phi(pr)\sim\phi(p)(\phi(q)+2)>\phi(pq)$. It implies that the difference between $\phi(qr)+\phi(pq)+\phi(q)$ and $\phi(qr)+\phi(pr)+\phi(r)$ is more than two. Likewise, same is true for $\phi(qr)+\phi(pq)+\phi(q)$ and $ \phi(qr)+\phi(pr)+\phi(pq).$ However, depending on the values of $p,q,r$, sometimes $\phi(qr)+\phi(pr)+\phi(r)$ is larger and sometimes $\phi(qr)+\phi(pr)+\phi(pq)$ is larger. Therefore, $x^{\phi(pr)+\phi(pq)+\phi(p)}+x^{\phi(qr)+\phi(pq)+\phi(q)}+x^{\phi(qr)+\phi(pr)+\phi(r)}$ or $x^{\phi(pr)+\phi(pq)+\phi(p)}+x^{\phi(qr)+\phi(pq)+\phi(q)}+x^{\phi(qr)+\phi(pr)+\phi(pq)}$ satisfies the property of log-concave. We are concerned with coefficients between $x^{\phi(qr)+\phi(pr)+\phi(r)} $ and $x^{\phi(qr)+\phi(pr)+\phi(pq)}.$ If $\phi(qr)+\phi(pr)+\phi(r)$ and $\phi(qr)+\phi(pr)+\phi(pq)$ differ by two, then
	 $$a_{\phi(qr)+\phi(pr)+\phi(r)+1}^{2}=0\ngeq a_{\phi(qr)+\phi(pr)+\phi(r)}a_{\phi(qr)+\phi(pr)+\phi(r)+2}=1,$$
	  where $\phi(qr)+\phi(pr)+\phi(pq)=\phi(qr)+\phi(pr)+\phi(r)+2$ or 
	  $$a_{\phi(qr)+\phi(pr)+\phi(pq)+1}^{2}=0\ngeq a_{\phi(qr)+\phi(pr)+\phi(pq)}a_{\phi(qr)+\phi(pr)+\phi(pq)+2}=1,$$ where $\phi(qr)+\phi(pr)+\phi(pq)+2=\phi(qr)+\phi(pr)+\phi(r).$ 
	  So, with the above observation, the polynomial given in \eqref{eq pqr} is log-concave if and only if the difference between $\phi(qr)+\phi(pr)+\phi(r) $ and $ \phi(qr)+\phi(pr)+\phi(pq)$ is away from two.
\end{proof}

As classified by the above result, we mention some of the classes of the independent domination polynomials of zero divisor graphs of $ \mathbb{Z}_{pqr}$ which are not log-concave.
\begin{enumerate}
	\item The independent domination polynomial of $G\cong \Gamma(\mathbb{Z}_{2\cdot 3\cdot 5})$ is
$ D_{i}(G,x)=x^{7}+x^{12}+x^{14}+x^{16}, $
which is not log-concave, since $a_{15}^{2}=0\ngeq a_{14}\cdot a_{16}=1.$	
\item The independent domination polynomial of $G\cong \Gamma(\mathbb{Z}_{2\cdot 5\cdot 7})$ is
$ D_{i}(G,x)=x^{11}+x^{32}+x^{34}+x^{36}, $
which is not log-concave, since $a_{35}^{2}=0\ngeq a_{34}\cdot a_{36}=1.$	
\item The independent domination polynomial of $G\cong \Gamma(\mathbb{Z}_{2\cdot 3\cdot 7})$ is
$ D_{i}(G,x)=x^{9}+x^{14}+x^{20}+x^{22}, $
which is not log-concave, since $a_{21}^{2}=0\ngeq a_{20}\cdot a_{22}=1.$	
\item The independent domination polynomial of $G\cong \Gamma(\mathbb{Z}_{3\cdot 5\cdot 7})$ is
$ D_{i}(G,x)=x^{22}+x^{36}+x^{42}+x^{44}, $
which is not log-concave, since $a_{43}^{2}=0\ngeq a_{42}\cdot a_{44}=1.$	
\end{enumerate}
For the above classes of polynomials, we see that $\phi(qr)+\phi(pr)+\phi(r) $ is bigger or  sometimes $ \phi(qr)+\phi(pr)+\phi(pq)$ is bigger, oscillation is exactly four, and they are not unimodal. 

Next, we discuss the independent domination polynomial of $\Gamma(\mathbb{Z}_{n})$ when $n$ is a prime power. We will state the result without proof, since it is similar to Theorem 7 given in \cite{gursoy}. We recall number theoretic facts: $\sum_{i=1}^{\ell}p^{i}=p^{\ell}-1$, where $p$ is prime and $\ell$ is a positive integer and for relatively primes $p$ and $q$, we have $\phi(pq)=\phi(p)\phi(q)$.
\begin{theorem}\label{thm zpalpha}
	Let $G\cong \Gamma(\mathbb{Z}_{n})$ be a zero divisor graph of $\mathbb{Z}_{n}$. Then the following hold.
	\begin{itemize}
		\item[\bf (i)] If $n=p^{2m}, $ where $p$ is prime and $m$ is a positive integer, then the independent domination polynomial of $G$ is given by
	\begin{align*}
		 D_{i}(G,x)&=\sum_{i=1}^{m} \phi(p^{i})\Big(x^{1+p^{2m-1}-p^{2m-1-(i-1)}}\Big)+x^{p^{2m-1}-p^{m}}.
	\end{align*}
	\item[\bf (ii)] If $n=p^{2m+1}, $  where $p$ is prime and $m$ is a positive integer, then the independent domination polynomial of $G$ is given by
	\begin{align*}
		D_{i}(G,x)&=\sum_{i=1}^{m} \phi(p^{i})\Big(x^{1+p^{2m}-p^{2m-(i-1)}}\Big)+x^{p^{2m}-p^{m}}.
	\end{align*}
	\end{itemize}
\end{theorem}

The following results gives the non unimodal and log-concave property of $\Gamma(\mathbb{Z}_{n})$ when $n$ is prime power.
\begin{theorem}\label{thm zpalpha restate}
	If $n$ is the prime power, then the independent domination polynomial of $\Gamma(\mathbb{Z}_{n})$ is log-concave but not unimodal. 
\end{theorem}
\noindent\begin{proof}
	Let $ n=p^{2m+1}$, where $p$ is prime and $m$ is a positive integer. Then the independent domination polynomial of $G\cong \Gamma(\mathbb{Z}_{n})$ is
	\begin{align*}
		D_{i}(G,x)&=\phi(p)x+\phi(p^{2})x^{1+\phi(p^{2m})}+\phi(p^{3})x^{1+\phi(p^{2m})+\phi(p^{2m-1})}+\dots+\\
		&\quad \phi(p^{m-1})x^{1+\phi(p^{2m})+\phi(p^{2m-1})+\dots+\phi(p^{m+3})}+\phi(p^{m})x^{1+\phi(p^{2m})+\phi(p^{2m-1})+\dots+\phi(p^{m+2})}\\
		&\quad+x^{\phi(p^{2m})+\phi(p^{2m-1})+\dots+\phi(p^{m+1})}.
	\end{align*}
	All the exponents of the above polynomial are odd except last exponent which is even. So, $2\leq \eta(D_{i}(G,x))$ and it follows that $D_{i}(G,x)$ is not unimodal. In fact $\eta(D_{i}(G,x))= m+1$ as all the coefficients of $x^{1+p^{2m}-p^{2m-(i-1)}}$ are strictly increasing, for $i=1,2,\dots,m$ and  $\phi(p^{2m})+\phi(p^{2m-1})+\dots+\phi(p^{m+1})=p^{2m}-p^{m}$ is strictly greater than $1+p^{2m}-p^{2m-(i-1)}$, for $i=1,2,\dots,m.$ So the coefficients of $D_{i}(G,x)$ increase precisely $m+1$ times.
	
	For $n=p^{2m},$ the independent domination polynomial of $G\cong \Gamma(\mathbb{Z}_{n})$ in standard form is
	\begin{align*}
		D_{i}(G,x)&=\phi(p)x+\phi(p^{2})x^{1+\phi(p^{2m-1})}+\phi(p^{3})x^{1+\phi(p^{2m-1})+\phi(p^{2m-2})}+\dots+\\
		&\quad \phi(p^{m-2})x^{1+\phi(p^{2m-1})+\phi(p^{2m-2})+\dots+\phi(p^{m+3})}+\phi(p^{m-1})x^{1+\phi(p^{2m-1})+\phi(p^{2m-2})+\dots+\phi(p^{m+2})}\\
		&\quad+x^{\phi(p^{2m-1})+\phi(p^{2m-2})+\dots+\phi(p^{m+1})}+\phi(p^{m})x^{1+\phi(p^{2m-1})+\phi(p^{2m-2})+\dots+\phi(p^{m+1})}.
	\end{align*}
	As in odd case, the exponents of the above expression are odd except $ \phi(p^{2m-1})+\phi(p^{2m-2})+\dots+\phi(p^{m+1})=p^{2m-1}-p^{m}$ and coefficients increase precisely $m+1$ times. From the representation of coefficients and exponents of $D_{i}(G,x)$,  it is easy to see that $1+p^{2m-1}-p^{2m-1-(i-1)}$  and $1+p^{2m}-p^{2m-(i-1)}$ have large gaps for respective values of $i.$ Thus by definition, the log-concave property holds trivially.
\end{proof}

\section{Conclusion}
\paragraph{}
This note gives the unimodal and log concave property of the independent domination polynomial of $\Gamma(\mathbb{Z}_{n})$ for $n\in \{p, p^{2}, 2p, pq, p^{2}q, pqr, p^{\alpha}\}$ and thereby corrects the results in \cite{gursoy} and flaws in their proof. The challenging task about the independent domination polynomial of $\Gamma(\mathbb{Z}_{n})$ is the location of zeros in the complex plane or possibly in some smaller annular region. The other graph polynomial like independent and the domination polynomials is another idea for carrying this work forward for the zero divisor graphs of commutative rings.

\section*{Data Availability:}
There is no data associated with this article.
\section*{Conflict of interest}
The authors declare that they have no competing interests.
\section*{Funding Statement}
There is no funding for this article.
\section*{Ethical Statement}
Not Applicable

\end{document}